\documentclass[11pt,a4paper]{article}
\usepackage[T1]{fontenc}
\usepackage[utf8]{inputenc}
\usepackage[margin=28mm]{geometry}
\usepackage{amsmath,amssymb,amsthm,mathtools, comment}
\usepackage{booktabs,array}
\usepackage[dvipsnames]{xcolor}
\usepackage[colorlinks=true,linkcolor=Blue,citecolor=Blue,urlcolor=Green]{hyperref}
\usepackage{tikz}
\usetikzlibrary{arrows.meta,decorations.pathreplacing}

\newtheorem{theorem}{Theorem}[section]
\newtheorem{lemma}[theorem]{Lemma}
\newtheorem{proposition}[theorem]{Proposition}
\newtheorem*{theorem*}{Theorem}
\newtheorem*{corollary*}{Corollary}

\theoremstyle{definition}
\newtheorem{definition}[theorem]{Definition}

\theoremstyle{remark}

\newcommand{\Z}{\mathbb Z}
\newcommand{\N}{\mathbb N}
\newcommand{\R}{\mathbb R}
\newcommand{\gp}[1]{\langle #1\rangle}
\newcommand{\dom}{\operatorname{dom}}
\newcommand{\rank}{\operatorname{rank}}
\newcommand{\Mod}{\operatorname{Mod}}
\newcommand{\Add}{\operatorname{Add}}
\newcommand{\Sub}{\operatorname{Sub}}

\newcommand{\PMod}{\operatorname{PMod}}
\newcommand{\Aut}{\operatorname{Aut}}

\newcommand{\calA}{\mathcal A}
\newcommand{\calR}{\mathcal R}

\newcommand{\calM}{\mathcal M}
\newcommand{\calD}{\mathcal D}
\newcommand{\calB}{\mathcal B}

\newcommand{\calC}{\mathcal C}

\newcommand{\doi}[1]{\href{https://doi.org/#1}{doi:\nolinkurl{#1}}}

\title{The membership problem for right-angled Artin groups and mapping class groups}
\author{
Ilya Kazachkov\\
{\small \textit{Ikerbasque - Basque Foundation for Science and Matematika Saila,}}\\
{\small \textit{UPV/EHU,  Sarriena s/n, 48940, Leioa - Bizkaia, Spain}}\\
{\small e-mail: \texttt{ilya.kazachkov@gmail.com}}
}
\begin{document}
\maketitle

\begin{abstract}
We prove that the membership problem in a right-angled Artin group $G(\Gamma)$ is decidable if and only if $G(\Gamma)$ is coherent, i.e. $\Gamma$ does not have induced $m$-cycles as subgraphs, for $m\ge 4$. More precisely, we show that any right-angled Artin group (raag) defined by an  $n$-cycle of length $n\ge 5$ contains a fixed finitely generated free subgroup with undecidable membership problem.

Consequently the membership problem in the mapping class group $M=\Mod(S)$ is decidable only if $M$ is virtually free or $S$ is the twice-punctured Klein bottle. In every other case, $M$ contains a fixed finitely generated free subgroup with undecidable membership problem.
\end{abstract}

\section{Introduction}
A landmark result in algorithmic group theory is Mihailova's construction of a finitely generated subgroup of $F_2\times F_2$ with undecidable membership problem \cite{Mihailova1958}. Mihailova's example demonstrates how difficult subgroup membership can be even when the ambient group has a simple presentation and a decidable word problem. Mihailova's example is finitely generated, but not finitely presented and any finitely presented subgroup of the direct product of $F_2\times F_2$ has decidable membership problem, see \cite[Theorem~2]{BaumslagRoseblade1984}.

In his paper \cite{Bridson2013}, Bridson proved that a right-angled Artin group can contain a finitely presented subgroup with undecidable membership problem \cite[Theorem~3.4]{Bridson2013}. He also established the corresponding phenomenon in mapping class groups of sufficiently large genus. Thus finite presentability of the subgroup does not suffice to ensure decidable membership. However, in Bridson's work both the ambient raag and the finitely presented subgroup are non-explicit, the raag contains $F_2 \times F_2$, and the construction does not yield subgroups with higher finiteness properties. (We note that direct products of free groups are particular examples of raags and, by the work of Bridson-Howie-Miller-Short, subgroups of direct products of free groups with higher finiteness properties enjoy much more rigid structure, see \cite[Theorem~A]{BHMS2002}.)

The above discussion raises two questions. Firstly, whether or not the undecidability of the membership problem for raags is rooted in the direct product of free groups and the Mihailova's construction and, secondly, how complex the subgroup with undecidable membership problem is, for instance in terms of its finiteness properties. In particular, Lohrey and Steinberg \cite[Section 5]{LohreySteinberg2008}, and more recently,  Gray and Nyberg-Brodda \cite[Section 4]{GrayNybergBrodda2025} ask if the membership problem for the raag defined by a pentagon is decidable.

In this paper, we characterise when a raag has (un)decidable membership problem and prove the following
\begin{theorem}[see Theorems \ref{thm:main} and \ref{thm:long-cycle-membership}]
Let $G(\Gamma)$ be a raag. Then the membership problem in $G(\Gamma)$ is decidable if and only if $\Gamma$ contains no induced cycles of length $\ge 4$ if and only if $G(\Gamma)$ is coherent.

Moreover, if $G(\Gamma)$ contains the raag $G(C_n)$ defined by the cycle of length $n\ge 5$, then it contains a fixed finitely generated free subgroup with undecidable membership problem.
\end{theorem}

The decidable direction is due to Kapovich, Weidmann and Myasnikov \cite[Corollary~1.3]{KapovichWeidmannMyasnikov2005}. In this paper we show that raags defined by cycles of length at least $5$ contain a finitely generated free group with undecidable membership problem.

We note that $G(\Gamma)$ contains the direct product of free groups if and only if $\Gamma$ contains an induced square. However, there are examples of graphs $\Gamma$ that do not contain induced cycles of length $\ge5$, but $G(\Gamma)$ contains the raag defined by the pentagon. For instance, raags defined by the complement graph of the cycle of length 6 or by the complement of a path on 8 vertices, see \cite[Theorem 4.2]{Kim2008} and \cite[Proposition 4.3]{CasalsRuizDuncanKazachkov2013}. 

\medskip

We then combine our construction for raags with known results on embeddings of raags into Artin groups of large type and mapping class groups and obtain a classification of such groups in terms of the (un)decidability of the membership problem.

\begin{theorem} \label{thm:orientable-surface-membership}
Let $S=S_{g,p}$ be a connected orientable surface of genus $g$, with $p$ punctures. Let $M=\Mod(S)$ be the orientation-preserving mapping class, with punctures allowed to permute. Set $\sigma(S)=3g-3+p$.
\begin{enumerate}
\item The group $M$ is virtually free if and only if $\sigma(S)\leq 1$. The subgroup membership problem in such $M$ is decidable.
\item If $\sigma(S)\geq 2$, then $M$ contains the raag $G$ defined by a pentagon as a subgroup. In particular, it contains a finitely generated free subgroup with undecidable membership problem.
\end{enumerate}
The surfaces with $\sigma(S)\leq 1$ are the spheres with $p\leq 4$ and the tori with $p\leq 1$.
\end{theorem}

\begin{theorem} \label{thm:nonorientable-surface-membership}
Let $N=N_{g,p}$ be a connected nonorientable surface of genus $g\ge1$, with $p\ge0$ punctures. Let $M=\Mod(N)$, with punctures allowed to permute.
\begin{enumerate}
\item
The subgroup membership problem in $M$ is decidable if and only if
$$
g=1,\ p\le3;
\qquad
g=2,\ p\le2;
\qquad\text{or}\qquad
(g,p)=(3,0).
$$
\item In every other case, $M$ contains the pentagon raag $G$. In particular, it contains a fixed finitely generated free subgroup with undecidable membership problem.
\end{enumerate}
\end{theorem}

The following corollary follows from \cite[Theorem 1]{CrispParis2001},  \cite[Theorem~1.1]{JankiewiczSchreve2022} (see also equation~(1.1)), \cite[Theorem 1.1]{Gordon2004}, \cite[Theorem~1]{FriedlWilton2016}, and Theorems \ref{thm:main} and \ref{thm:long-cycle-membership}

\begin{corollary*}
Let $A_\Gamma$ be an Artin group. 
\begin{enumerate}
\item
If the $2$-subgraph of $\Gamma$ contains an induced cycle of length at least $5$, then $A_\Gamma$ contains a fixed finitely generated free subgroup with undecidable membership problem.
\item
If $A_\Gamma$ is of large type,  then the membership problem in $A_\Gamma$ is decidable if and only if $\Gamma$ contains no cycles of length $\ge 3$.  If  $\Gamma$ contains a cycle of length $\ge 3$, then $A_\Gamma$ contains a fixed finitely generated free subgroup with undecidable membership problem.
\end{enumerate}
\end{corollary*}

\subsection*{Outline of the proof}

The proof has two main ingredients. We begin in Section~\ref{sec:two-counter-machines} by fixing a deterministic two-counter machine for which acceptance is undecidable, following Sapir \cite{Sapir2015}. The first ingredient is the construction in Section~\ref{sec:progression-graphs}, where we encode any deterministic two-counter machine by finitely many pairs of arithmetic progressions.

More precisely, any set of pairs of arithmetic progressions defines a graph whose vertex set is $\Z$ and two vertices $n$ and $m$ are joined by an edge if and only if there is an integer $S$ and a pair of arithmetic progressions $\{\alpha+pk, k\in \Z\}$ and $\{\beta+qk, k\in \Z\}$ in the set so that $n=\alpha+pS$ and $m=\beta+qS$. This arithmetic construction is related to the works of Conway and Leary, \cite{Conway1972, Leary2021}.

We prove that the machine accepts an input precisely when the integer encoding its initial configuration lies in the connected component of the integer encoding its accepting configuration. Applying this construction to the fixed machine gives a fixed graph and a fixed vertex $t$ for which no algorithm decides whether or not $n$ lies in the connected component of $t$.

The second ingredient is encoding the pairs of arithmetic progressions into the decidability of an algebraic problem for the raag $P$ defined by the path with $4$ vertices $a, b,c, d$.

A crucial easy observation is that in a cylinder covering a torus (with fixed horizontal and vertical generators) a horizontal path and a vertical path that start at two different marked points meet at the same vertex precisely when the horizontal and vertical displacements form two arithmetic progressions whose terms are naturally paired.

We want to transfer the pairing of arithmetic progressions in the cylinder to connectivity by paths in a special square complex. In order to do so, we take one cylinder for each pair of arithmetic progressions, one auxiliary cylinder, and two planes $\pi_1$ and $\pi_2$. We then identify the horizontal and vertical lines in the cylinders with some lines in the planes. This gives us a special square complex and its fundamental group $V$ is a finitely generated free subgroup of $P$. With the prescribed identification of lines, we show that a horizontal path $b^n$ in the square complex belongs to the double coset $V\cdot \langle a, d \rangle$ in $P$ if and only if $n$ lies in the connected component of $0$ in the graph defined by pairs of arithmetic progressions. Hence we obtain that the membership in  double cosets of subgroups in $P$ is undecidable:
\begin{theorem} \label{thm:path-raag-product-membership}
Let $P=\gp{a,b,c,d\mid [a,b]=[b,c]=[c,d]}$ be the raag defined by the path on $4$ vertices. Then there exists a fixed finitely generated free subgroup $V\le P$ such that there is no algorithm to decide for $n\in\Z$ whether or not $b^{n}\in V\cdot \gp{a,d}$.

Thus membership in the product of the two fixed finitely generated free subgroups $V$ and $\gp{a,d}$ is undecidable, even when the input is restricted to powers of $b$.
\end{theorem}

We then use a trick with centralizer extensions to promote the undecidability of the membership for double cosets in $P$ to the undecidability of the membership problem in the pentagon raag. Write the raag defined by the pentagon as follows
$$
G=\gp{P,e\mid e^{-1}\ell e=\ell,\quad\ell\in \gp{a,d}}.
$$
Thus $G$ is an HNN-extension of $P$ with a stable letter that centralizes $L=\gp{a,d}$. Then, as we prove in Section \ref{sec:centralizer-extension}, $geg^{-1}\in\gp{V,e}$  if and only if $g\in VL$ which is undecidable by Theorem \ref{thm:path-raag-product-membership}.

\begin{theorem}\label{thm:main}
There exists a fixed finitely generated free subgroup $H\le G$ so that there is no algorithm to decide whether or not $b^{-n}eb^n$ belongs to $H$, hence $G$ has undecidable membership problem.
\end{theorem}

Finally, in Section~\ref{sec:longer-cycles} we extend the construction to raags defined by cycles of length greater than five. In Section~\ref{sec:mapping-class-groups} we apply our results to mapping class groups.

\subsection*{Acknowledgment}
This work was supported by the Basque Government grant IT1913-26, and PID2025-171325NB-I00 and PCI2024-155053-2 grants funded by MICIU/AEI/10.13039/501100011033 and FEDER/EU.

I am grateful to Montserrat Casals-Ruiz for the feedback on earlier versions of the manuscript.

\subsection*{AI disclosure statement}
OpenAI's ChatGPT was used to draw figures and proof verification.  Proofs of Theorems \ref{thm:orientable-surface-membership} and \ref{thm:nonorientable-surface-membership} boil down to references and were generated by AI. The author assumes full responsibility for the contents of this publication.

\section{Universal two-counter machine}
\label{sec:two-counter-machines}

We follow \cite[Sections 2.1 and 2.2]{Sapir2015} and will use the universality theorem, \cite[Theorem~2.3(a)]{Sapir2015}. 

Minsky or two-counter machine has finitely many states and two counters with values in $\N$. Its states are labelled $0,1,\ldots,M-1$, with start state $1$ and distinguished stop state $0$. A configuration is a triple $(s,x_1,x_2)$. The counters are numbered $1,2$, containing $x_1,x_2$, respectively. Put $\mathbf x=(x_1,x_2)$, $\mathbf e_1=(1,0)$, and $\mathbf e_2=(0,1)$. 

A program is a finite collection of commands, each acting as a partially defined map on the configuration space. The domain of a command consists of the configurations at its source state that satisfy its stated condition, if any. A command with source state $s\ne 0$ has one of the following forms, see \cite[Section 2.2]{Sapir2015}, where $c\in\{1,2\}$ and $r\in\{0,\ldots,M-1\}$:
$$
\begin{aligned}
\Add(c);r:\quad
&(s,\mathbf x)\longmapsto(r,\mathbf x+\mathbf e_c),\\
x_c>0\to\Sub(c);r:\quad
&(s,\mathbf x)\longmapsto(r,\mathbf x-\mathbf e_c)
\quad(x_c>0),\\
x_c=0\to r:\quad
&(s,\mathbf x)\longmapsto(r,\mathbf x)
\quad(x_c=0).
\end{aligned}
$$

The other counter is unchanged. The last two commands are defined only when their indicated conditions hold. There are no commands with source state $0$.

A machine is \emph{deterministic} if distinct commands have disjoint domains as partial maps on configurations (that is, there is no ambiguity in what machine does). Their union is a partial \emph{step function}
\begin{equation}\label{eq:step_function}
\delta_{\calM}: \{0,\ldots,M-1\}\times\N^2 \rightharpoonup
\{0,\ldots,M-1\}\times\N^2.
\end{equation}
The domain of $\delta_{\calM}$ consists exactly of those configurations at which a command is applicable. A configuration with non-zero state may lie outside this domain. A run is obtained by iterating $\delta_{\calM}$ for as long as it is defined. We distinguish \emph{acceptance}, meaning that the run reaches $(0,0,0)$, from stopping because no command applies at some other configuration. Thus each configuration has at most one successor.

\begin{theorem}[{Two-counter universality; \cite[Theorem~2.3(a)]{Sapir2015}}]
\label{thm:counter-halting}
For every recursively enumerable set $X\subseteq\N$, there is a fixed deterministic two-counter machine $\calM$ such that, for every $m\in\N$, its run from $(1,2^m,0)$ reaches $(0,0,0)$ if and only if $m\in X$. If $m\notin X$, that run is infinite.
\end{theorem}

In particular, the specified input runs never stop at a  configuration other than $(0,0,0)$.

\section{Arithmetic progressions and two-counter machines}\label{sec:progression-graphs}

In this section, we translate acceptance by the fixed two-counter
machine into connectivity in a graph (whose vertices are labelled by $\Z$). We encode configurations of the machine as integers and represent the machine's transitions by finitely many arithmetic-progression rules with pairwise disjoint source domains. Lemma \ref{lem:sink} then shows that the integer that encodes an input of the machine is connected to the integer that encodes its terminal configuration if and only if the machine accepts that input. This defines a fixed graph specified by finitely many rules, in which membership in the connected component of a fixed vertex is undecidable.

This arithmetic construction is related to the works of Conway and Leary. Conway \cite{Conway1972} showed that iteration of simple arithmetic rules can simulate arbitrary computation. Leary \cite[Theorem~11]{Leary2021} used graphs built from the iterates of a function to prove that connectedness is undecidable for recursively described graphs. In our setting, the question is whether a given integer belongs to one specified component of a single graph.

\begin{definition}\label{def:rules}
A \emph{progression rule} is a quadruple
$$
\rho=(\alpha,p,\beta,q), \qquad \alpha,\beta\in\Z,\quad 0<p,q\in \N.
$$
A set of progression rules defines a non-oriented graph $\calA(\calR)$ as follows. The vertex set of $\calA(\calR)$ is $\Z$. To each progression rule we associate a set of edges
$$
E_\rho= \left\{
(\alpha+pk,\beta+qk), k\in\Z,\ \alpha+pk\ne\beta+qk
\right\}.
$$
Thus, a finite collection $\calR$ of progression rules, defines the undirected graph
$$
\calA(\calR)= \left(\Z,\ \bigcup_{\rho\in\calR} E_\rho\right).
$$
Write $n\sim_{\calR}m$ when $n$ and $m$ belong to the same connected component of this graph. In the definition of $E_\rho$ we explicitly omit loops (note that loops do not affect connectivity).

Given a progression rule $\rho=(\alpha,p,\beta,q)$, we refer to the sets $\alpha+p\Z$ and $\beta+q\Z$ as source and target domains or progressions.
\end{definition}

A progression rule $\rho=(\alpha,p,\beta,q)$ defines a (partial) map
$$
T_\rho:\alpha+p\Z\longrightarrow\beta+q\Z,
\qquad
T_\rho(x)=\beta+\frac qp(x-\alpha).
$$
The domain of $T_\rho$ is exactly $\alpha+p\Z$. For $x\notin\alpha+p\Z$, the value $T_\rho(x)$ is undefined, even if $\beta+\frac qp(x-\alpha)$ is an integer. 

Any partial function $f:X\rightharpoonup X$ defines a directed graph $\calD(f)$ as follows. The vertex set of $\calD(f)$ is $X$. For every $x\in\dom(f)$ in the domain of $f$, there is a directed edge $x\longrightarrow f(x)$. Let $\calA(f)$ be the undirected graph obtained from $\calD(f)$ by forgetting the directions of edges.

For each $x\in\dom(T_\rho)$, the rule $\rho$ maps $x$ to $T_\rho(x)$.
We call this step from $x$ to $T_\rho(x)$ a \emph{directed transition} of $\rho$. For a finite collection of rules $\calR$, forgetting the directions of all the transitions, discarding loops and repeated edges gives the graph $\calA(\calR)$.

In a connected directed graph a vertex that has no outgoing edges
is called a \emph{sink}. In $\calD(f)$ a vertex $t$ is a sink if $t\notin\dom(f)$. The following lemma describes the connected component of a sink in $\calA(f)$. It follows immediately from \cite[Lemma~2.1]{Sapir2015} which says that  two configurations of a deterministic machine are equivalent under forward and backward commands precisely when they both admit forward computations ending in the same configuration. Since $f$ is undefined at $t$, the only defined forward iterate of $t$ is $f^0(t)=t$. Thus $v$ and $t$ have a common forward iterate if and only if $f^k(v)=t$ for some $k\ge0$. Formally,

\begin{lemma}[Sink lemma]\label{lem:sink}
Let $f:X\rightharpoonup X$ be a partial function and suppose
$t\in X\setminus\dom(f)$. 

Then the connected component of $t$ in $\calA(f)$ is
$$
\calB_t=\{x\in X: f^j(x)=t\text{ for some }j\ge0\}.
$$
Here $f^0$ is the identity map, and the iterates of $f$ must be defined.
\end{lemma}

\begin{proposition}\label{prop:arithmetic}
For every deterministic two-counter machine $\mathcal M$, one can effectively construct a finite set $\calR=\{\rho_1,\ldots,\rho_r\}$ of progression rules and a computable injective map $E$ from the configurations of $\mathcal M$ to $\Z$ with the following properties.

The source domains $\alpha_i+p_i\Z$ of the rules in $\calR$ are pairwise disjoint, so their directed transitions define a partial function $f:\Z\rightharpoonup\Z$. For every configuration $\xi$ and every terminal configuration $\eta$ of $\mathcal M$, we have
$$
E(\xi)\sim_{\calR}E(\eta)
\quad\Longleftrightarrow\quad
\text{$\mathcal M$ reaches $\eta$ from $\xi$ after finitely many steps}.
$$
Here zero steps are allowed, and terminal means that no command applies to $\eta$. Moreover, $E(\eta)\notin\dom(f)$.

In particular, by {\rm Theorem~\ref{thm:counter-halting}}, every nonrecursive recursively enumerable set gives rise to a fixed finite set $\calR$ of progression rules and a fixed integer $t$ such that there is no algorithm to decide for $n\in\Z$ whether or not $n\sim_{\calR}t$. The source domains $\alpha_i+p_i\Z$ of the rules in $\calR$ are pairwise disjoint, and $t\notin\dom(f)$ for the corresponding partial function $f$.
\end{proposition}
\begin{proof}
Fix a deterministic two-counter machine $\calM$. Let $M$ be its number of states, with start state $1$ and stop state $h=0$.

We encode machine configurations by the injective map
$$
E:\{0,\ldots,M-1\}\times\N^2\longrightarrow\N,
\qquad
E(s,u,v)=M2^u3^v+s.
$$

The rule list $\calR$ is determined by the fixed machine.  Each command contributes one or two progression rules representing its action on encoded configurations.

Suppose a command at state $s$ is $\Add(1);s'$. Then the machine transitions as follows $(s,u,v)\longmapsto(s',u+1,v)$. On the level of encodings, we have:
$$
E(s,u,v)=s+M2^u3^v \qquad \text{and} \qquad
E(s',u+1,v)=s'+2M2^u3^v.
$$
The command at $s$ determines the progression rule $\rho=(s,M,s',2M)$. Recall that $\rho$  defines a (partial) map $T_\rho$.
For a configuration $(s,u,v)$, taking $k=2^u3^v$ gives
$$
T_\rho(E(s,u,v))=E(s',u+1,v) \qquad T_\rho(s+Mk)=s'+2Mk, \quad k\in\Z
$$
And so $T_\rho$ reproduces the corresponding command. Note that the domain of $T_\rho$ is $s+M\Z$, which also contains integers that do not encode configurations.

Similarly, assign to each command with source state $s$ the  progression rules specified in the following table. Here $s'$ is its successor state, $u,v$ are the values of counters $1,2$, and $k$ ranges over $\Z$.
\begin{center}
\renewcommand{\arraystretch}{1.35}
\begin{tabular}{@{}ll@{}}
\toprule
Command at state $s$ & Directed transition\\
\midrule
$\Add(1);s'$ & $s+Mk\longrightarrow s'+2Mk$\\
$u>0\to\Sub(1);s'$ & $s+2Mk\longrightarrow s'+Mk$\\
$u=0\to s'$ & $M+s+2Mk\longrightarrow M+s'+2Mk$\\
$\Add(2);s'$ & $s+Mk\longrightarrow s'+3Mk$\\
$v>0\to\Sub(2);s'$ & $s+3Mk\longrightarrow s'+Mk$\\
$v=0\to s'$ &
$M\rho+s+3Mk\longrightarrow M\rho+s'+3Mk$,\\
& $\rho\in\{1,2\}$\\
\bottomrule
\end{tabular}
\end{center}
The row for $v=0$ specifies two rules, one for each $\rho\in\{1,2\}$. No rule is assigned to the stop state $h$. Let $\calR$ be the resulting finite rule list.

We first show that the source domains of the progression rules in $\calR$ are pairwise disjoint. Every rule arising at state $s$ has source domain contained in $s+M\Z$, so rules arising at different states have disjoint source domains.

At a fixed state, an increment command excludes every other command and contributes just one progression rule. Otherwise, determinism permits at most two conditional commands, testing the same counter: one when it is zero and the other when it is positive. Indeed, any other pair of conditional commands would apply simultaneously to some configuration.

For counter $1$, the possible source progressions are
\begin{equation} \label{eq:first-counter-source-progressions}
s+2M\Z,
\qquad
s+M+2M\Z.
\end{equation}
For counter $2$, they are
\begin{equation} \label{eq:second-counter-source-progressions}
s+3M\Z,
\qquad
s+M+3M\Z,
\qquad
s+2M+3M\Z,
\end{equation} 
where the last two arise from a single zero-test command. The progressions in \eqref{eq:first-counter-source-progressions} are pairwise disjoint, and so are the progressions in \eqref{eq:second-counter-source-progressions}. Since conditional commands at a fixed state test the same counter, all source domains are pairwise disjoint.

Set $f(x)=T_\rho(x)$ for $x\in\dom(T_\rho)$. This defines a partial function $f:\Z\rightharpoonup\Z$ with
$$
\dom(f)
=\bigcup_{\rho\in\calR}\dom(T_\rho)
\subseteq\Z\setminus M\Z.
$$
In particular, $M\notin\dom(f)$.

We now verify that $f$ reproduces the machine's transitions on encoded configurations and is undefined there exactly when no command applies.

For a configuration $\xi=(s,u,v)$, set $N=2^u3^v$. The conditions on the counters can be expressed arithmetically as
$$
\begin{aligned}
u=0 &\iff 2\nmid N,
&\qquad u>0 &\iff 2\mid N,\\
v=0 &\iff 3\nmid N,
&\qquad v>0 &\iff 3\mid N.
\end{aligned}
$$
The definitions of the directed transitions, see the table, and the above equivalences show that, for each command with source state $s$, the integer $E(\xi)$ belongs to the source domain of one of its associated progression rules if and only if that command applies to $\xi$.

Recall that $\delta_{\calM}$ is the machine's partial step function, see Equation \eqref{eq:step_function}. Whenever a command applies, the formulas in the table imply
$$
f(s+MN)=s'+MN',
$$
where $N'=2N$ or $N'=3N$ for an increment of counter $1$ or $2$, respectively; $N'=N/2$ or $N'=N/3$ for a decrement of the  corresponding counter; and $N'=N$ for a zero-test. Thus $f(E(\xi))$ is the encoding of the configuration obtained by applying the command to $\xi$. If no command applies, $E(\xi)$ belongs to none of the source progressions, whether $s=0$ or the machine is blocked at a nonzero state. Therefore, for every configuration $\xi$,
$E(\xi)\in\dom(f)$ if and only if $\xi\in\dom(\delta_{\calM})$, and, whenever these conditions hold, we have  $f(E(\xi))=E(\delta_{\calM}(\xi))$.

By induction on $j$ we see that $f^j(E(\xi))$ is defined if and only if so is $\delta_{\calM}^j(\xi)$, and 
\begin{equation}\label{eq:iterates}
f^j(E(\xi))=E(\delta_{\calM}^j(\xi)).
\end{equation}

By construction, $\calA(\calR)$ is the undirected graph of $f$, with loops and repeated edges omitted. These omissions do not affect connected components.

Let $\xi$ be any configuration and let $\eta$ be a terminal configuration. We have shown above that $E(\eta)\notin\dom(f)$. So,  by Lemma~\ref{lem:sink}, $E(\xi)\sim_{\calR}E(\eta)$ if and only if $f^j(E(\xi))=E(\eta)$, for some $j\ge 0$. By \eqref{eq:iterates}, the latter happens if and only $\calM$ reaches $\eta$ from $\xi$.

Finally, fix a nonrecursive recursively enumerable set $X$ and apply the construction to the machine $\calM$ provided by Theorem~\ref{thm:counter-halting}. For input $m\in\N$, its initial configuration is $\xi_0=(1,2^m,0)$ and its accepting configuration is $(0,0,0)$. Their encodings are
$$
n(m):=E(1,2^m,0)=M2^{\,2^m}+1,
\qquad
t:=E(0,0,0)=M.
$$
By Theorem~\ref{thm:counter-halting}, it follows that $m\in X$ if and only if  $n(m)\sim_{\calR}t$.

The set $\calR$ and the integer $t$ are fixed. Deciding membership in the connected component of $t$ implies $X$ is recursive.
\end{proof}

\section{Realization by double cosets}\label{sec:double-coset-realization}

We give a brief outline of this section. In Section~\ref{sec:cylinder-progressions} we associate a cylinder with each pair of arithmetic progressions. The period fixes their differences, and the position of the vertical path fixes the offsets. In Section \ref{sec:construction-of-Y} we assemble these cylinders into a square complex $Y$ mapping to the standard square complex of $P$. All cylinders attach to the same two square-tiled planes, namely the $(b,a)$-plane $\Pi_A$ and the $(d,c)$-plane $\Pi_D$. Each cylinder is attached along a chosen $b$-labelled line to $\Pi_A$ and along a chosen $c$-labelled line to $\Pi_D$. Paths using only the labels $a$ and $d$ pass between these planes through the intersections defined by the progression rules.

There are two points to establish. First, $Y$ must define a subgroup of $P$ with the claimed free basis. We do this in Sections \ref{sec:free-basis-of-u} and \ref{sec:local-isometry}. Second, the full covering associated with that subgroup must not introduce additional $a,d$-connections between the marked vertices, see Section \ref{sec:ad-paths-and-double-cosets}.

\subsection{Arithmetic progressions in cylinders}
\label{sec:cylinder-progressions}
We start with an easy observation on how to encode a pair of arithmetic progressions in a cylinder. This serves as a warm-up and to establish notation.

Let $T$ be the Salvetti complex of $\langle b,c\mid[b,c]=1\rangle\cong\Z^2$. The universal cover is $\R^2$, where $b$ and $c$ label the horizontal and vertical edges, respectively.

Fix positive integers $p,q$. The quotient by the translation $(p,-q)$ is a square complex homeomorphic to a cylinder. We write $|x,y|$ for the image of the point $(x,y)$ in this quotient. The 
covering map to $T$ is then given by
$$
\calC_{p,q}=\R^2/\langle(p,-q)\rangle
\longrightarrow T,
\qquad
|x,y|\longmapsto(x,y)+\Z^2.
$$
Then $|x,y|=|x',y'|$ if and only if $(x'-x,y'-y)=k(p,-q)$ for some $k\in\Z$.

This is the covering corresponding to the cyclic subgroup
$\langle h\rangle$, where $h=b^p c^{-q}$. The lift of $h$ starting at $o=|0,0|$ is a closed path going once around the cylinder. Hence, $h$ records the period $(p,-q)$.

The paths labelled $b^n$ and $c^m$ starting at $o$ end at the same vertex precisely when $|n,0|=|0,m|$ or whenever $(n,-m)=k(p,-q)$ for some $k\in\Z$. Thus the possible pairs of exponents are exactly $(n,m)=(pk,qk)$, with $k\in\Z$. The direction of the period determines the ratio $q/p$, while the period vector determines the increments $p$ and $q$. Observe that the possible pairs of exponents define a pair of arithmetic progressions $pk$ and $qk$, $k\in \Z$.

We now allow the vertical path to start at a marked vertex.  Fix integers $\alpha,\beta$, and let
$$
o'=|\alpha,-\beta|\in \calC{p,q},\qquad \Delta=b^\alpha c^{-\beta}.
$$
The word $\Delta$ labels a fixed path from $o$ to $o'$. For integers $n,m$, let $u(n)$ be the endpoint of the horizontal path $b^n$ starting at $o$, and let $v(m)$ be the endpoint of the vertical path $c^m$ starting at $o'$.  In coordinates,
$$
u(n)=|n,0|,\qquad v(m)=|\alpha,m-\beta|.
$$

The endpoints of these paths coincide exactly when the two displayed lifts differ by a whole number of cylinder periods. Subtracting their coordinates we have 
\begin{equation}  \label{eq:cylinder-progressions}
 u(n)=v(m)
 \quad\Longleftrightarrow\quad
 (n-\alpha,\beta-m)=k(p,-q)  
 \quad\Longleftrightarrow\quad  n=\alpha+pk,\quad m=\beta+qk.
\end{equation}

Summarising, the cylinder $\calC{p,q}$ with a marked point $o'=|\alpha,-\beta|$ uniquely determines a pair of arithmetic progressions  $\{\alpha+pk,k\in\Z\}$ and $\{\beta+qk, k\in \Z\}$. The period of $\calC{p,q}$ records the differences and the element $\Delta$ records the offsets of the progressions.

The terms of the corresponding progressions are naturally paired by \eqref{eq:cylinder-progressions}, $\{\, \alpha+pk,\beta+qk):k\in\Z\,\}$.

We note that we do not assume that $p$ and $q$ are co-prime. If $r=\gcd(p,q)>1$, the vector $(p/r,-q/r)$ has the same direction as $(p,-q)$, but it is not the period of the corresponding cylinder.

\subsection{Example}

We begin with a single pair of arithmetic progressions encoded in a cylinder. Using this cylinder, we build a square complex that translates connectivity in the progression graph into a double-coset membership problem in $P$. In the next section, we extend this construction to any finite family of pairs.

We illustrate the construction with the rule $\rho=(1,2,0,1)$. We start with two square-tiled planes
$$
\Pi_A=\mathbb R^2_{b,a},
\qquad
\Pi_D=\mathbb R^2_{d,c}.
$$
Here $b,d$ and $a,c$ label the corresponding horizontal and vertical edges.

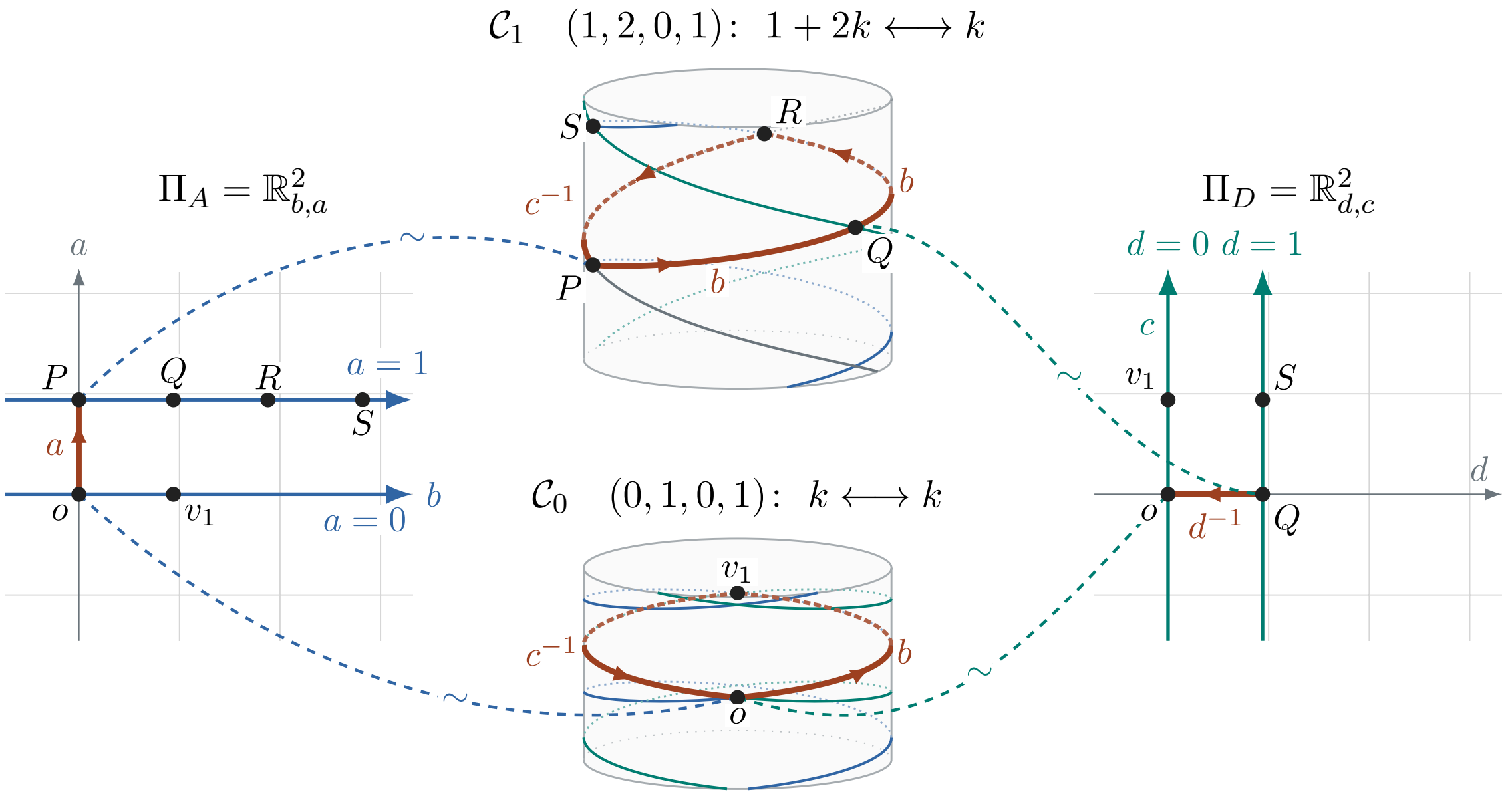
\begin{figure}[ht]
\centering
\begin{tikzpicture}[x=1cm,y=1cm,>=Latex,font=\small,
  ygrid/.style={draw=black!16,line width=.36pt},
  yB/.style={draw=yblue,line width=1.0pt},
  yC/.style={draw=yteal,line width=1.0pt},
  yglue/.style={dash pattern=on 2.5pt off 2.4pt,line width=.85pt},
  yarrow/.style={-{Latex[length=2.2mm,width=1.6mm]}},
  ylab/.style={fill=white,inner sep=1pt,font=\small},
  ypoint/.style={circle,fill=black!87,inner sep=1.5pt}]
\definecolor{yblue}{RGB}{48,101,164}
\definecolor{yteal}{RGB}{0,125,113}
\definecolor{ywalk}{RGB}{157,64,32}
\definecolor{yslate}{RGB}{107,117,124}
\def\Ymark#1#2#3{\node[ypoint] at (#1) {};
  \node[ylab,outer sep=2pt,#2] at (#1) {$#3$};}
\def\Ystrokes#1#2#3#4#5#6{%
  \pgfmathtruncatemacro{\Ylo}{floor((#1)/180)}
  \pgfmathtruncatemacro{\Yhi}{floor((#2)/180)}
  \foreach \Yhalf in {\Ylo,...,\Yhi}{
    \pgfmathsetmacro{\Yu}{max(#1,180*\Yhalf)}
    \pgfmathsetmacro{\Yv}{min(#2,180*(\Yhalf+1))}
    \ifdim\Yv pt>\Yu pt
      \pgfmathtruncatemacro{\Yback}{mod(abs(\Yhalf),2)}
      \ifnum\Yback=0 \tikzset{thisY/.style={#6}}\else
                       \tikzset{thisY/.style={#5}}\fi
      \draw[thisY, line join=round] plot[domain=\Yu:\Yv,samples=25,variable=\ang]
       ({\yrad*cos(\ang)},{(#4)*(\ang-(#3))/360+\yellipse*sin(\ang)});
    \fi}}
\def\Ytrack#1#2#3{%
  \pgfmathsetmacro{\Ystart}{(#1)+360*\ybottom/(#2)}
  \pgfmathsetmacro{\Ystop}{(#1)+360*\ytop/(#2)}
  \pgfmathsetmacro{\Yfirst}{min(\Ystart,\Ystop)}
  \pgfmathsetmacro{\Ylast}{max(\Ystart,\Ystop)}
  \Ystrokes{\Yfirst}{\Ylast}{#1}{#2}
    {draw=#3,line width=.75pt}
    {draw=#3!55,dash pattern=on .6pt off 1.1pt,line width=.56pt}}
\def\Yedge#1#2#3#4#5{%
  \Ystrokes{#1}{#2}{#3}{#4}
    {draw=ywalk,line width=1.65pt}
    {draw=ywalk!83,dash pattern=on 2pt off 1pt,line width=1.15pt}
  \pgfmathsetmacro{\Ymid}{(#1)+(#5)*((#2)-(#1))}
  \draw[draw=ywalk,line width=1.0pt,yarrow]
    plot[domain={\Ymid-7}:\Ymid,samples=5,variable=\ang]
    ({\yrad*cos(\ang)},{(#4)*(\ang-(#3))/360+\yellipse*sin(\ang)});}
\def\Yshell{%
  \fill[black!2] (-\yrad,\ybottom)--(-\yrad,\ytop)
    arc[start angle=180,end angle=0,x radius=\yrad,y radius=\yellipse]
    --(\yrad,\ybottom)
    arc[start angle=0,end angle=-180,x radius=\yrad,y radius=\yellipse]--cycle;
  \draw[yslate!60,line width=.55pt] (-\yrad,\ybottom)--(-\yrad,\ytop)
    (\yrad,\ybottom)--(\yrad,\ytop)
    (0,\ytop) ellipse[x radius=\yrad,y radius=\yellipse]
    (-\yrad,\ybottom) arc[start angle=180,end angle=360,x radius=\yrad,y radius=\yellipse];
  \draw[yslate!45,dotted] (-\yrad,\ybottom)
    arc[start angle=180,end angle=0,x radius=\yrad,y radius=\yellipse];}
\def\Yplanegrid{%
  \begin{scope}
  \clip (-.78,-1.55) rectangle (3.53,2.35);
  \draw[ygrid] (-2,-2) grid (4,3);
  \end{scope}}
\begin{scope}[shift={(1.00,3.05)},x=.94cm,y=.94cm]
  \Yplanegrid
  \draw[yB,->] (-.78,0)--(3.53,0);
  \draw[yB,->] (-.78,1)--(3.53,1);
  \draw[yslate,->] (0,-1.55)--(0,2.39) node[above] {$a$};
  \node[text=yblue,right] at (3.53,0) {$b$};
  \node[text=yblue,ylab,below left] at (3.51,-0.1) {$a=0$};
  \node[text=yblue,ylab,above=6pt] at (3.27,1) {$a=1$};
  \coordinate (Ao) at (0,0); \coordinate (AP) at (0,1);
  \coordinate (AQ) at (1,1); \coordinate (AR) at (2,1);
  \coordinate (Av) at (1,0); \coordinate (AS) at (3,1);
  \draw[ywalk,line width=1.65pt] (Ao)--(AP);
  \draw[ywalk,line width=1pt,yarrow] (0,.40)--(0,.78);
  \node[ylab,text=ywalk,left=3pt] at (0,.5) {$a$};
\end{scope}
\begin{scope}[shift={(11.83,3.05)},x=.94cm,y=.94cm]
  \Yplanegrid
  \draw[yC,->] (0,-1.55)--(0,2.39);
  \draw[yC,->] (1,-1.55)--(1,2.39);
  \draw[yslate,->] (-.78,0)--(3.53,0) node[above left] {$d$};
  \node[above,text=yteal] at (0,2.39) {$d=0$};
  \node[above,text=yteal] at (1,2.39) {$d=1$};
  \node[text=yteal,left] at (0,1.77) {$c$};
  \coordinate (Do) at (0,0); \coordinate (DQ) at (1,0);
  \coordinate (Dv) at (0,1); \coordinate (DS) at (1,1);
  \draw[ywalk,line width=1.65pt] (DQ)--(Do);
  \draw[ywalk,line width=1pt,yarrow] (.70,0)--(.32,0);
  \node[ylab,below=3pt,text=ywalk] at (.51,0) {$d^{-1}$};
\end{scope}
\node[font=\normalsize] at (2.64,6.04) {$\Pi_A=\mathbb R^2_{b,a}$};
\node[font=\normalsize] at (13.03,6.04) {$\Pi_D=\mathbb R^2_{d,c}$};
\def\yrad{1.53}\def\yellipse{.29}
\begin{scope}[shift={(7.55,5.43)}]
  \def\ybottom{-1.04}\def\ytop{1.56}
  \Yshell
  \Ytrack{200}{1.38}{yblue}
  \Ytrack{380}{-2.76}{yteal}
  \Ytrack{200}{-2.76}{yslate}
  \coordinate (CP) at ({\yrad*cos(200)},{\yellipse*sin(200)});
  \coordinate (CQ) at ({\yrad*cos(320)},{.46+\yellipse*sin(320)});
  \coordinate (CR) at ({\yrad*cos(440)},{.92+\yellipse*sin(440)});
  \coordinate (CS) at ({\yrad*cos(560)},{1.38+\yellipse*sin(560)});
  \Yedge{200}{320}{200}{1.38}{.40}
  \Yedge{320}{440}{200}{1.38}{.78}
  \Yedge{440}{560}{560}{-2.76}{.44}
  \node[ylab,text=ywalk,below] at (-.20,-.09) {$b$};
  \node[ylab,text=ywalk,right] at (1.56,.75) {$b$};
  \node[ylab,text=ywalk,left] at (-1.55,.55) {$c^{-1}$};
\end{scope}
\node[font=\normalsize] at (7.55,7.69) {$\mathcal C_1\quad (1,2,0,1)\colon\ 1+2k\longleftrightarrow k$};
\begin{scope}[shift={(7.55,1.32)}]
  \def\ybottom{-.91}\def\ytop{1.00}
  \Yshell
  \Ytrack{270}{.92}{yblue}
  \Ytrack{270}{-.92}{yteal}
  \coordinate (Co) at (0,-\yellipse);
  \coordinate (Cv) at (0,{.46+\yellipse});
  \Yedge{270}{450}{270}{.92}{.34}
  \Yedge{450}{630}{630}{-.92}{.77}
  \node[ylab,text=ywalk,right] at (1.54,.17) {$b$};
  \node[ylab,text=ywalk,left] at (-1.54,.17) {$c^{-1}$};
\end{scope}
\node[font=\normalsize] at (7.55,3.02) {$\mathcal C_0\quad (0,1,0,1)\colon\ k\longleftrightarrow k$};
\draw[yglue,yblue] (AP) .. controls (2.44,5.48) and (4.66,5.97) .. (CP) node[pos=.64,ylab,text=yblue] {$\sim$};
\draw[yglue,yteal] (CQ) .. controls (10.00,5.88) and (11.04,3.17) .. (DQ) node[pos=.57,ylab,text=yteal] {$\sim$};
\draw[yglue,yblue] (Ao) .. controls (2.9,1.29) and (5.28,.47) .. (Co) node[pos=.59,ylab,text=yblue] {$\sim$};
\draw[yglue,yteal] (Co) .. controls (9.6,.40) and (10.39,1.55) .. (Do) node[pos=.51,ylab,text=yteal] {$\sim$};
\Ymark{Ao}{below left}{o}
\Ymark{Av}{below right}{v_1}
\Ymark{AP}{above left}{P}
\Ymark{AQ}{above}{Q}
\Ymark{AR}{above}{R}
\Ymark{AS}{below}{S}
\Ymark{Do}{below left}{o}
\Ymark{Dv}{above left}{v_1}
\Ymark{DQ}{below right}{Q}
\Ymark{DS}{above right}{S}
\Ymark{CP}{below left}{P}
\Ymark{CQ}{below right}{Q}
\Ymark{CR}{above right}{R}
\Ymark{CS}{left}{S}
\Ymark{Co}{below}{o}
\Ymark{Cv}{above}{v_1}
\end{tikzpicture}
\caption{The square complex $Y_\rho$ for the rule $\rho=(1,2,0,1)$. The upper cylinder pairs $1+2k$ with $k$, and the lower, auxiliary cylinder pairs $k$ with itself. The brown edges trace a free basis of $\pi_1(Y,o) =\gp{z_0=bc^{-1}, z_1=a(b^2c^{-1})a^{-1}, w_1=abd^{-1}}$}
\label{fig:one-rule-cylinder-attachments}
\end{figure}

We first consider the auxiliary cylinder identifying the vertices on the $b$- and $c$- axis. Formally, consider the rule $\rho_0=(0,1,0,1)$ pairing each integer $k$ with itself. By the preceding section, it determines a cylinder $\calC_0=\R^2/\langle(1,-1)\rangle$. We view $\calC_0$ as quotient of the square-tiled plane $\R^2_{b,c}$ ($b$ corresponding to the horizontal axis and $c$ vertical). Attach the $b$-line through $|0,0|$ in $\calC_0$ to the horizontal line defined by $a=0$ in $\Pi_A$ and the $c$-line through $|0,0|$ to the vertical line defined by $d=0$ in $\Pi_D$, see Figure \ref{fig:one-rule-cylinder-attachments}.  For each $n\in\Z$, we have $|n,0|=|0,n|$ in $\mathcal C_0$. Thus $(n,0)\in\Pi_A$ and $(0,n)\in\Pi_D$ give the same vertex after gluing. Denote this vertex by $v_n$ and take $o=v_0$ as the basepoint.

Now consider the rule $\rho=(1,2,0,1)$, which pairs $1+2k$ with $k$. As before, $\rho$ determines a cylinder $\calC_1=\R^2_{b,c}/\langle(2,-1)\rangle$ with marked point $o_1'=|1,0|=Q$.  We attach the $b$-line in $\calC_1$ through $P$ to the horizontal line $a=1$ in $\Pi_A$, and the  $c$-line through $Q$ to the vertical line $d=1$ in $\Pi_D$, see Figure \ref{fig:one-rule-cylinder-attachments}.  In particular, the $b$-edge from $P$ to $Q$ joins the two attaching lines in $\calC_1$.

Let $Y_\rho$ be the obtained square complex. Then the fundamental group $\pi_1(Y_\rho,o)$ is free with the following set of generators
\begin{itemize}
\item $z_0=bc^{-1}$ winds around the lower cylinder $\calC_0$;
\item $z_1=a(b^2c^{-1})a^{-1}$ goes from $o$ to $P$, winds around the upper cylinder $\calC_1$ and goes back to $o$;
\item $w_1=abd^{-1}$ follows the path $o\to P\to Q\to o$ through $\Pi_A$, $\calC_1$ and back through $\Pi_D$.
\end{itemize}
The edge labels define a combinatorial local isometry from $Y_\rho$ to the Salvetti complex of the path raag $P$. The induced homomorphism from $\pi_1(Y_\rho,o)$ to $P$ is therefore injective.

To see how the rule appears in $Y_\rho$, fix $n\in\mathbb Z$. The rule pairs $2n+1$ with $n$ and $|2n+1,0|=|1,n|$ in $\calC_1$, by \eqref{eq:cylinder-progressions}. Hence $(2n+1,1)\in\Pi_A$ and $(1,n)\in\Pi_D$ are the same vertex of $Y_\rho$. Starting at $v_n$
in $\Pi_D$, we can therefore follow the path
$$
v_n=(0,n)_{\Pi_D}
\xrightarrow{\ d\ }
(1,n)_{\Pi_D}=(2n+1,1)_{\Pi_A}
\xrightarrow{\ a^{-1}\ }
(2n+1,0)_{\Pi_A}=v_{2n+1}.
$$
Thus $da^{-1}$ takes $v_n$ to $v_{2n+1}$. The path labelled by $(da^{-1})^2$ connects $v_n$ to $v_{2n+1}$ to $v_{2(2n+1)+1}$, and so on.

\subsection{Construction of \texorpdfstring{$Y$}{Y}}\label{sec:construction-of-Y}
We first construct the square complex $Y$. Fix a finite list of progression rules
$$
\calR=\{(\alpha_i,p_i,\beta_i,q_i),i=1, \ldots, r\}
\qquad
\alpha_i,\beta_i\in\Z,\quad 0<p_i,q_i\in\N.
$$
For each $1\le i\le r$, take a separate copy of the marked cylinder from Section~\ref{sec:cylinder-progressions},
$$
\calC_i=\R^2_{b,c}/\Z(p_i,-q_i),
\qquad
o_i=|0,0|,\qquad o_i'=|\alpha_i,-\beta_i|.
$$
Here $|x,y|$ denotes the image of the point $(x,y)\in\mathbb R^2_{b,c}$ in $\mathcal C_i$. The $b$-axis is horizontal and the $c$-axis is vertical.

For the construction of $Y$, also set
$$
(\alpha_0,p_0,\beta_0,q_0)=(0,1,0,1).
$$
The corresponding auxiliary cylinder is
$$
\mathcal C_0=\mathbb R^2_{b,c}/\mathbb Z(1,-1),
\qquad o_0=o_0'=|0,0|.
$$
The auxiliary rule pairs each integer with itself and hence does not change the relation $\sim_{\mathcal R}$.

All the cylinders attach to the same two square-tiled planes, $\Pi_A=\mathbb R^2_{b,a}$, and $\Pi_D=\mathbb R^2_{d,c}$. The subscripts give the horizontal and vertical directions, in that
order. Denote the origins of $\Pi_A$ and $\Pi_D$ by $o_A$ and $o_D$,
respectively.

For each $0\le i\le r$, let $\ell_b^i$ be the $b$-labelled line through $o_i$ in $\mathcal C_i$, and let $\ell_c^i$ be the $c$-labelled line through $o_i'$. We attach $\ell_b^i$ to the horizontal line $a=i$ in $\Pi_A$ and $\ell_c^i$ to the vertical line $d=i$ in $\Pi_D$ identifying
\begin{equation}\label{eq:cylinder-attachments}
|n,0|\in\mathcal C_i
\ \sim\ (n,i)\in\Pi_A,\quad \text{and}\quad
|\alpha_i,m-\beta_i|\in\mathcal C_i
\ \sim\ (i,m)\in\Pi_D.
\end{equation}
Here $n,m\in\mathbb R$ parametrize the two lines, and integer values give their vertices. Taking $n=m=0$, we see that $o_i$ is attached to $(0,i)\in\Pi_A$ and $o_i'$ to $(i,0)\in\Pi_D$. Different cylinders attach to distinct lines in each plane. Let $Y$ be the resulting square complex.

In the auxiliary cylinder $\mathcal C_0$, the equality $|s,0|=|0,t|$ holds precisely when $(s,-t)$ is an integer multiple of $(1,-1)$, that is, when $s=t\in\mathbb Z$. Thus, for each integer $n$, the auxiliary cylinder identifies the vertex $(n,0)\in\Pi_A$ with the vertex $(0,n)\in\Pi_D$. Denote their common image in $Y$ by $v_n$. In particular, the origins $o_A$ and $o_D$ give the same vertex $o=v_0$, which we take as the basepoint of $Y$.

\begin{figure}[ht]
\centering
\begin{tikzpicture}[
  x=.91cm,y=.91cm,
  font=\small,
  >=Latex,
  gridline/.style={draw=black!13,line width=.3pt},
  axis/.style={draw=black!45,line width=.55pt,->},
  bline/.style={draw=blue!70!black,line width=1.25pt},
  cline/.style={draw=teal!80!black,line width=1.25pt},
  route/.style={draw=black,line width=1.6pt,-{Latex[length=2.1mm]}},
  glue/.style={densely dashed,line width=.9pt},
  vertex/.style={circle,fill=black,inner sep=1.8pt},
  lab/.style={fill=white,inner sep=2pt}
]
\begin{scope}
  \clip (.35,.05) rectangle (4.15,3.65);
  \foreach \x in {.8,1.35,1.9,2.45,3,3.55,4.1}
    \draw[gridline] (\x,.05)--(\x,3.65);
  \foreach \y in {.2,.7,1.2,1.7,2.2,2.7,3.2}
    \draw[gridline] (.35,\y)--(4.15,\y);
\end{scope}
\begin{scope}
  \clip (10.85,.05) rectangle (14.95,3.65);
  \foreach \x in {11.2,11.7,12.2,12.7,13.2,13.7,14.2,14.7}
    \draw[gridline] (\x,.05)--(\x,3.65);
  \foreach \y in {.2,.7,1.2,1.7,2.2,2.7,3.2}
    \draw[gridline] (10.85,\y)--(14.95,\y);
\end{scope}

\node[font=\normalsize] at (2.25,4.55)
  {$\Pi_A=\mathbb R^2_{b,a}$};
\node[font=\normalsize] at (7.5,4.55) {$\mathcal C_i$};
\node[font=\normalsize] at (12.9,4.55)
  {$\Pi_D=\mathbb R^2_{d,c}$};

\draw[axis] (.35,.7)--(4.15,.7) node[right] {$b$};
\draw[axis] (.8,.05)--(.8,3.8) node[above] {$a$};
\draw[bline,->] (.35,2.7)--(4.1,2.7);
\node[text=blue!70!black,lab,above=4pt] at (1.85,2.7) {$a=i$};
\node[vertex,fill=blue!70!black] at (.8,2.7) {};
\node[lab,above left=3pt] at (.8,2.7) {$(0,i)$};
\draw[route] (3,.7)--(3,2.7)
  node[midway,lab,left=4pt] {$a^i$};
\node[vertex] at (3,.7) {};
\node[vertex] at (3,2.7) {};
\node[lab,below=5pt] at (3,.7) {$(n,0)$};
\node[lab,above=5pt] at (3,2.7) {$(n,i)$};

\fill[black!2] (6.15,0)--(6.15,3.2)
  arc[start angle=180,end angle=0,x radius=1.35,y radius=.3]
  --(8.85,0)
  arc[start angle=0,end angle=-180,x radius=1.35,y radius=.3];
\draw[black!35] (6.15,0)--(6.15,3.2);
\draw[black!35] (8.85,0)--(8.85,3.2);
\draw[black!35] (7.5,3.2) ellipse[x radius=1.35,y radius=.3];
\draw[black!35] (6.15,0)
  arc[start angle=180,end angle=360,x radius=1.35,y radius=.3];
\draw[black!20,dotted] (6.15,0)
  arc[start angle=180,end angle=0,x radius=1.35,y radius=.3];
\node[align=center,font=\footnotesize] at (7.5,3.93)
  {local view of the\\two attaching lines};
\draw[bline,->] (6.45,1.2)--(8.55,1.2);
\draw[cline,->] (7.5,.3)--(7.5,2.8);
\node[text=blue!70!black,lab,above=3pt] at (8.35,1.2) {$b$};
\node[text=teal!80!black,lab,right=3pt] at (7.5,2.65) {$c$};
\node[vertex] at (7.5,1.2) {};
\node[lab,above left=2pt] at (7.5,1.2) {$\nu$};

\draw[axis] (10.85,.7)--(15,.7) node[right] {$d$};
\draw[axis] (12.2,.05)--(12.2,3.8) node[above] {$c$};
\draw[cline,->] (14.2,.05)--(14.2,3.8);
\node[text=teal!80!black,above] at (14.2,3.8) {$d=i$};
\node[vertex,fill=teal!80!black] at (14.2,.7) {};
\node[lab,below right=4pt] at (14.2,.7) {$(i,0)$};
\draw[route] (14.2,2.7)--(12.2,2.7)
  node[midway,lab,above=4pt] {$d^{-i}$};
\node[vertex] at (12.2,2.7) {};
\node[vertex] at (14.2,2.7) {};
\node[lab,above right=4pt] at (14.2,2.7) {$(i,m)$};
\node[lab,below left=4pt] at (12.2,2.7) {$(0,m)$};

\draw[glue,blue!70!black]
  (3,2.7) .. controls (4.6,1.1) and (6,-.55) .. (7.5,1.2)
  node[pos=.5,lab] {$\sim$};
\draw[glue,teal!80!black]
  (7.5,1.2) .. controls (9,-.55) and (13.2,1.1) .. (14.2,2.7)
  node[pos=.5,lab] {$\sim$};
\node[vertex] at (3,2.7) {};
\node[vertex] at (7.5,1.2) {};
\node[vertex] at (14.2,2.7) {};

\node[align=center] at (7.5,-.9)
  {$\nu=|n,0|=|\alpha_i,m-\beta_i|$\\[4pt]
   $n=\alpha_i+p_i k,\qquad m=\beta_i+q_i k$};
\end{tikzpicture}
\caption{The attachment of cylinders $\mathcal C_i$. Suppose $n=\alpha_i+p_i k$ and $m=\beta_i+q_i k$. The two cylinder endpoints are then the same vertex $\nu$.  After gluing, the paths labelled $a^i$ and $d^{-i}$ concatenate.
}
\label{fig:cylinder-attachments}
\end{figure}
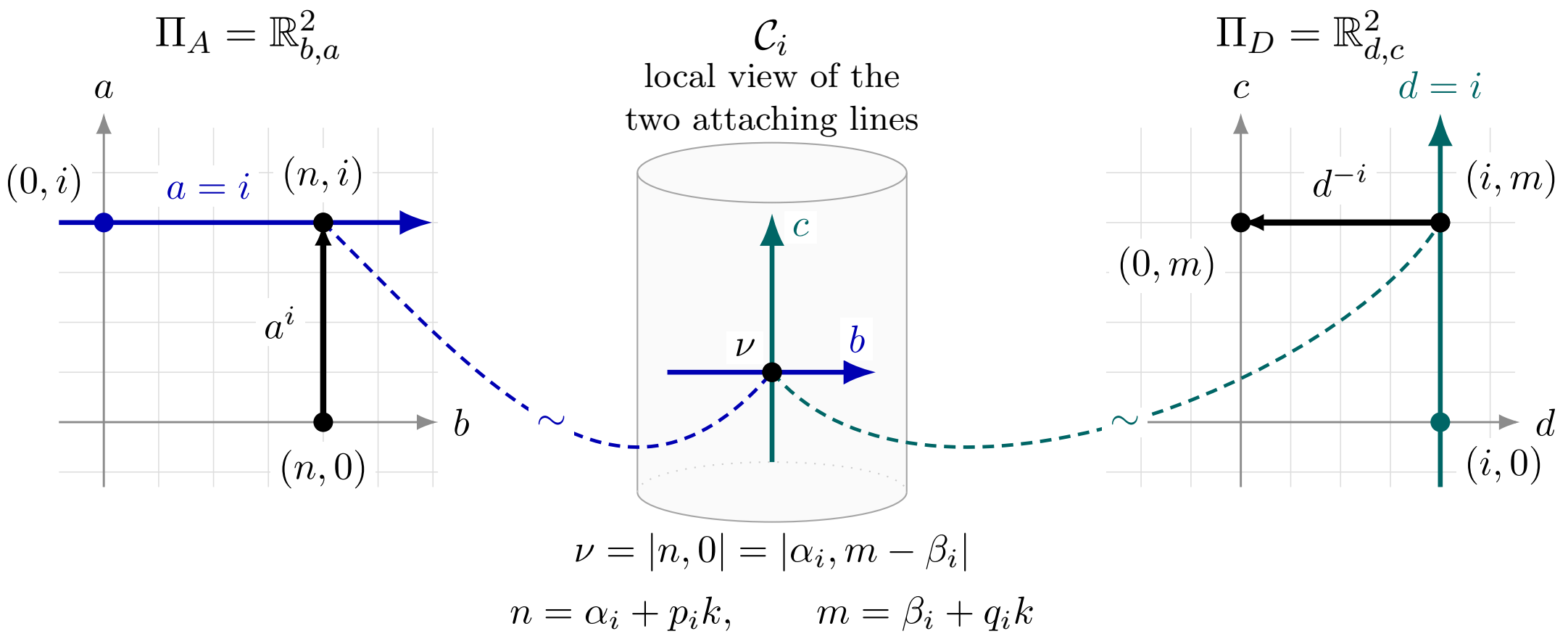

Now suppose that the $i$-th rule pairs $n$ with $m$, where $1\le i\le r$. Thus $n=\alpha_i+p_i k$ and $m=\beta_i+q_i k$
for some $k\in\mathbb Z$. By \eqref{eq:cylinder-progressions}, we have that $|n,0|=|\alpha_i,m-\beta_i|$ in $\calC_i$. It follows from \eqref{eq:cylinder-attachments} that $(n,i)\in\Pi_A$ and $(i,m)\in\Pi_D$ are the same vertex of $Y$. At this vertex the following two paths concatenate, see Figure~\ref{fig:cylinder-attachments}:
$$
v_n=(n,0)_{\Pi_A}
\xrightarrow{\ a^i\ }
(n,i)_{\Pi_A}=(i,m)_{\Pi_D}
\xrightarrow{\ d^{-i}\ }
(0,m)_{\Pi_D}=v_m.
$$
Thus, a single application of the $i$-th rule gives a path from $v_n$ to $v_m$ labelled $a^i d^{-i}$. Following this path in reverse gives a path from $v_m$ to $v_n$ labelled $d^i a^{-i}$. Since the auxiliary cylinder identifies the two vertices representing each integer, the paths for successive rule applications concatenate in $Y$.

The attachments preserve edge labels and orientations, so the covering maps from the planes and cylinders to their respective tori define a map $Y\to X$, where $X$ is the Salvetti complex of $P$.

Just as in \eqref{eq:cylinder-progressions}, the two attaching lines in $\mathcal C_i$ meet at the parameters $n,m$ precisely when
\begin{equation}\label{eq:exact-transition}
\begin{aligned}
|n,0|=|\alpha_i,m-\beta_i|
&\quad\Longleftrightarrow\quad
(n-\alpha_i,\beta_i-m)=k(p_i,-q_i)
\text{ for some }k\in\Z\\
&\quad\Longleftrightarrow\quad
n=\alpha_i+p_i k,\qquad m=\beta_i+q_i k
\text{ for some }k\in\Z.
\end{aligned}
\end{equation}

We now formulate the main result of this section. Define 
$$
h_i=b^{p_i}c^{-q_i},  \Delta_i=b^{\alpha_i}c^{-\beta_i}, \qquad  
z_i=a^ih_i a^{-i}, w_j=a^j\Delta_j d^{-j} \qquad (0\le i\le r,1\le j\le r).
$$
Their roles are as follows.
\begin{itemize}
\item The word $h_i$ records the period $(p_i,-q_i)$ of the
$i$-th cylinder. Its lift goes once around that cylinder, exactly
as the word $h$ in Section~\ref{sec:cylinder-progressions}.
\item The word $\Delta_i$ records the offset path from $|0,0|$ to
$|\alpha_i,-\beta_i|$. Its endpoint is the start of the vertical
path used to compare the two progressions. This is the word $\Delta$
from Section~\ref{sec:cylinder-progressions}, with the parameters of the $i$-th
rule; its lift need not be closed.
\item The word $z_i$ gives a loop at the common basepoint of $Y$:
follow $a^i$ to the $i$-th cylinder, go around it along $h_i$, and
return along $a^{-i}$.
\item The word $w_i$ gives a second kind of based loop: follow
$a^i$ to the cylinder, follow $\Delta_i$ to its attachment to the right
plane, and return to the basepoint along $d^{-i}$.
\end{itemize}
Thus $h_i,\Delta_i$ describe paths within a cylinder, while $z_i,w_i$
are labels of based loops that generate the fundamental group. For the auxiliary cylinder, $h_0=z_0=bc^{-1}$ and $\Delta_0=1$, its two attaching lines meet at the common basepoint.

\begin{proposition}[Arithmetic realization]\label{prop:realization}
The group
$$
U=\langle z_0,\ldots,z_r,w_1,\ldots,w_r\rangle\le P
$$
is free on the indicated generators, $U\cong F_{2r+1}$.

For all $n,m\in\Z$, we have that $n\sim_{\calR}m$ if and only if $Ub^nL=Ub^mL$.
\end{proposition}
We prove this proposition in Section \ref{sec:ad-paths-and-double-cosets}.

\subsection{\texorpdfstring{$\pi_1(Y)$}{The fundamental group of Y} is free} \label{sec:free-basis-of-u}

In this section we show that $Y$ is homotopically equivalent to a graph and its fundamental group has rank $2r+1$. 

First attach all the cylinders to $\Pi_A$. Each cylinder has infinite cyclic fundamental group and is attached along a contractible line. Hence, this defines a free group  of rank $r+1$. We now attach $\Pi_D$. Its intersection with the complex already constructed consists of $r+1$ pairwise disjoint lines. Since $\Pi_D$ is contractible, this attachment adds $r$ free generators and we conclude that $\pi_1(Y,o)\cong F_{2r+1}$.

The free basis can be read directly from the construction. To go once around $\mathcal C_i$, first follow $a^i$ from $o$ to $|0,0|$, then follow $b^{p_i}c^{-q_i}$ around the cylinder and return along $a^{-i}$. This loop has label $z_i$. For each $1\le i\le r$, another loop follows $a^i$ from $o$ to $|0,0|$, then follows $\Delta_i=b^{\alpha_i}c^{-\beta_i}$ to the right attaching line, and returns to $o$ along $d^{-i}$ in $\Pi_D$. Its label is $w_i$. Using $\mathcal C_0$ to connect the two planes gives the free basis 
$\{z_0,\ldots,z_r,w_1,\ldots,w_r\}$.
The labels of these loops generate the image $U$ of $\pi_1(Y,o)$ in
$P$. We next show that the map of the fundamental groups is injective.

\subsection{Local isometry of \texorpdfstring{$Y$ to $X$}{Y to X}}
\label{sec:local-isometry}

The vertex link in $X$ is flag, so $X$ is nonpositively curved.  It suffices to check that $Y\to X$ is a local isometry. For a square complex mapping to $X$, this means that the map on each vertex link is injective and its image is a full subgraph of the link in $X$. Equivalently, distinct incident edges have distinct labels, and any pair of edges labelled by commuting generators spans a square. This follows directly by construction.

By \cite[Lemma~2.11]{HaglundWise2008}, the local isometry induces an embedding $\pi_1(Y,o)\hookrightarrow\pi_1(X)=P$. On the level of universal covers, we obtain an embedding of $\widetilde Y$ as a convex subcomplex of $\widetilde X$.

Let $X_U\to X$ be the covering corresponding to $U$. Taking the quotient of the equivariant embedding $\widetilde Y  \hookrightarrow \widetilde X$ by $U$, we obtain an embedding $Y\hookrightarrow X_U$ inducing the isomorphism $\pi_1(Y,o)\cong U=\pi_1(X_U,U)$. Vertices of $X_U$ are the right cosets $Ug$, and an edge labelled $s$ connects $Ug$ to $Ugs$. Under this identification, $o$ is $U$ and $v_n$ is $Ub^n=Uc^n$, the paths starting at $o$ and labelled $b^n$ in $\Pi_A$ and $c^n$ in $\Pi_D$ end at the same vertex of the auxiliary cylinder.

\subsection{Putting everything together}
\label{sec:ad-paths-and-double-cosets}

Call an edge path in $Y$ an $(a,d)$-path if its label is a word in 
$a^{\pm1},d^{\pm1}$. 
\begin{lemma}\label{lem:ad-connectivity}
The vertices $v_n$ and $v_m$ are joined by an $(a,d)$-path in $Y$
if and only if  $n\sim_{\calR}m$.
\end{lemma}
\begin{proof}
By Section~\ref{sec:construction-of-Y}, an application of rule $\rho_i$ defines a path labelled $a^id^{-i}$ between the corresponding marked vertices. Concatenating these paths proves $\Leftarrow$.

Conversely, write $\lambda_n$ for the $a$-line with $b$-coordinate $n$ in $\Pi_A$, and $\mu_m$ for the $d$-line with $c$-coordinate $m$
in $\Pi_D$. An $(a,d)$-path can pass from one line to another only
where they meet. At $v_j$ it can pass between $\lambda_j$ and $\mu_j$. At the intersection in any other cylinder, the two indices are paired by that cylinder's rule, see \eqref{eq:exact-transition}. Thus, along a path from $v_n$ to $v_m$, each change of index is an application of a rule, hence the lemma.
\end{proof}

\begin{proof}[Proof of Proposition~\ref{prop:realization}]
A freely reduced $(a,d)$-path in $X_U$ with endpoints in $Y$ lies in $Y$. Indeed, its lift to $\widetilde X$ is geodesic and the preimage of $Y$ is convex by Section~\ref{sec:local-isometry}. It therefore contains the lifted path.

Since $v_n=Ub^n$, an $(a,d)$-path in $X_U$ joins $v_n$ to $v_m$ precisely when $Ub^n\ell=Ub^m$ for some $\ell\in L$. Hence, by Lemma \ref{lem:ad-connectivity}, we have that $Ub^nL=Ub^mL$ if and only if  $n\sim_{\calR}m$.

Sections~\ref{sec:free-basis-of-u} and~\ref{sec:local-isometry} describe the free basis of $U$.
\end{proof}

\begin{proof}[Proof of Theorem~\ref{thm:path-raag-product-membership}]
Choose $\calR,t$ as in Proposition~\ref{prop:arithmetic} and let $V=b^{-t}Ub^t$. Then, by definition of $V$, $b^{-k}\in VL$  if and only if  $Ub^{t-k}L=Ub^tL$. By Proposition~\ref{prop:realization}, this is equivalent to $t-k\sim_{\calR}t$. The subgroup $V$ is fixed and free of finite rank. Deciding the first condition is equivalent to deciding membership in the component of $t$.
\end{proof}

\section{Transfer to the centralizer extension}\label{sec:centralizer-extension}

The following normal-form argument applies to arbitrary base groups. We use the HNN normal-form theorem and refer to  \cite[Chapter~IV, Section~2]{LyndonSchupp2001}.

\begin{proposition}\label{prop:transfer}
Let $K$ be a group, $S\le K$, and
$$
\widehat K=\gp{K,e\mid e^{-1}se=s\quad(s\in S)}.
$$
For every subgroup $V\le K$ and every $g\in K$,
$$
geg^{-1}\in\gp{V,e}
\quad\Longleftrightarrow\quad
g\in VS.
$$
Note that we do not assume that $V\cap S$ is finitely generated.
\end{proposition}
\begin{proof}
Put $J=V\cap S$. Let 
$$
\widehat{V}=\gp{V,\tau\mid \tau^{-1}j\tau=j\quad(j\in J)}.
$$
The map $\tau\longmapsto e$ induces a homomorphism from $\widehat V$ to $\widehat K$. This homomorphism is injective on $V$ by definition. Furthermore, it is injective on the whole group by Britton's lemma. Indeed, a reduced word in $\widehat{V}$ remains reduced in $\widehat K$, since a coefficient in $V$ belongs to $S$ if and only if it belongs to $J$. The image of $\widehat{V}$ is $\gp{V,e}$.

Take $g=vs$, where $v\in V$ and $s\in S$, then
$$
geg^{-1}=vse s^{-1}v^{-1}=vev^{-1}\in\gp{V,e}.
$$

Conversely, suppose that $geg^{-1}\in\gp{V,e}$.
This element has HNN-length one in $\widehat K$. Since reduced words in $\widehat{V}$ remain reduced in $\widehat K$, by normal forms in HNN extensions, we have that $geg^{-1}=v_0ev_1$, where $v_0,v_1\in V$.

Rearranging we conclude that $e^{-1}(v_0^{-1}g)e=v_1g\in K$. By Britton's lemma, $v_0^{-1}g\in S$, and therefore $g\in VS$.
\end{proof}

\section{Membership in \texorpdfstring{$G$}{G}}
\label{sec:pentagon-membership}

\begin{proposition}[Freeness at a terminal configuration]
\label{prop:terminal-freeness}
Suppose that the source progressions of the rules in $\calR$ are
pairwise disjoint and that $t$ lies in none of them. Let
$U_t=b^{-t}Ub^t$ and $H=\langle U_t,e\rangle$. Then
$U_t\cap L=\{1\}$.

Consequently, $H\cong U_t*\langle e\rangle$ is free of rank $2r+2$.
A free basis consists of $e$, the elements $b^{-t}z_i b^t$
for $0\le i\le r$, and the elements $b^{-t}w_i b^t$ for
$1\le i\le r$.
\end{proposition}
\begin{proof}
We first show that $U_t\cap L=\{1\}$. Recall that $\lambda_n$ is the $a$-line with $b$-coordinate $n$ in $\Pi_A$, and $\mu_n$ is the $d$-line with $c$-coordinate $n$ in $\Pi_D$. The auxiliary cylinder $\mathcal C_0$ identifies $(n,0)\in\lambda_n$ with $(0,n)\in\mu_n$. Their common image is $v_n$, the unique vertex of $\mathcal C_0$ on $\lambda_n$ or $\mu_n$. We call the vertices $v_n$ \emph{auxiliary vertices}.

Any other intersection of $\lambda_n$ with a $d$-line comes from one of the original rules. More precisely, by \eqref{eq:exact-transition}, the $i$-th cylinder identifies $(n,i)\in\lambda_n$ with $(i,m)\in\mu_m$ precisely when $n=\alpha_i+p_i k$ and $m=\beta_i+q_i k$ for some $k\in\Z$. We call this vertex a \emph{rule intersection}. Since the source progressions are pairwise disjoint, at most one rule applies to $n$, and for that rule $k$ is uniquely defined. Thus $\lambda_n$ contains at most one rule intersection. In particular, the only vertex of $\lambda_t$ at which a path can switch to a $d$-line is $v_t$.

Suppose that $1\ne u\in U_t\cap L$. Then, the freely reduced word $u=u(a,d)$ labels a closed path at $Ub^t=v_t$ in $X_U$, and this path lies in $Y$ by Section~\ref{sec:ad-paths-and-double-cosets}. Divide it into maximal $a$- and $d$- segments labelled by nonzero powers of $a$ and $d$ correspondingly. Then, the first segment is a $d$-segment: a nonzero $a$-segment starting at $v_t$ moves away from $v_t$ along $\lambda_t$. Therefore, it can not switch to a $d$-line or return to $v_t$. Applying the same argument to the reverse path shows that the last segment is also a $d$-segment.

We claim that every $d$-segment of this path starts at an auxiliary
vertex. The first starts at $v_t$. Suppose a $d$-segment starts at
$v_m$ and is followed by an $a$-segment along $\lambda_n$. The switch
takes place in a vertex distinct from $v_m$, so it is a rule intersection. It is therefore the unique rule intersection on $\lambda_n$. The $a$-segment moves away from this point, so its only possible switch back to a $d$-line is at $v_n$. This proves the claim.

Therefore, the last $d$-segment must start at an auxiliary vertex and end at the auxiliary vertex $v_t$. This is impossible, since any $d$-line contains only one auxiliary vertex. Hence $U_t\cap L=\{1\}$.

Conjugating the basis of $U$ from Proposition~\ref{prop:realization} by $b^{-t}$ and adjoining $e$ gives the free basis of rank $2r+2$.
\end{proof}

\begin{proof}[Proof of Theorem~\ref{thm:main}]
Take the same $\calR,t$ and $U_t=b^{-t}Ub^t$ as in the proof of
Theorem~\ref{thm:path-raag-product-membership}. By Proposition~\ref{prop:terminal-freeness}, $H=\langle U_t,e\rangle$ is a fixed finitely generated free subgroup of $G$. By Proposition~\ref{prop:transfer}, for every $k\in\Z$ one has that $b^{-k}eb^k\in H$ if and only if $b^{-k}\in U_tL$. The latter is undecidable by Theorem~\ref{thm:path-raag-product-membership}.
\end{proof}

\section{RAAGs defined by longer cycles}\label{sec:longer-cycles}

Fix $N\ge 6$ and let $s=N-5$. Label the vertices of the $N$-cycle in
cyclic order by
$$
a,\ b,\ c_0,\ c_1,\ldots,c_s,\ d,\ e.
$$
Let $G_N$ be its raag and $P_N$ be the subgroup generated by all vertices except $e$. Then $P_N$ is the raag defined by the path from $a$ to $d$. The only neighbours of $e$ are $a$ and $d$. Thus, with $L=\langle a,d\rangle$,
$$
 G_N=\langle P_N,e\mid e^{-1}\ell e=\ell\quad(\ell\in L)\rangle
     =P_N*_L\bigl(L\times\langle e\rangle\bigr).
$$
In Section~\ref{sec:double-coset-realization} the rule $\rho_i$ defines a single cylinder which is attached onto the planes $\Pi_A$ and $\Pi_D$. Here we use a chain of $s+1$ cylinders for each
$i$. The first encodes the same rule, and all the others just transfer the information: they identify the vertex with integer parameter $m$ on its incoming line with the vertex having the same parameter on its outgoing line. 

\subsection{The square complex \texorpdfstring{$Y_N$}{Y\_N} for the longer path}

Fix a list $\mathcal R$ of $r$ rules as in Section~\ref{sec:construction-of-Y}. Use $(\alpha_0,p_0,\beta_0,q_0)=(0,1,0,1)$ for the auxiliary rule. Take the endpoint planes
$\Pi_A=\mathbb R^2_{b,a}$ and $\Pi_D=\mathbb R^2_{d,c_s}$. For each
$0\le i\le r$, form a separate chain of cylinders
$$
 \mathcal C_{i,0}=\mathbb R^2_{b,c_0}/\mathbb Z(p_i,-q_i),
 \qquad
 \mathcal C_{i,j}=\mathbb R^2_{c_{j-1},c_j}/\mathbb Z(1,-1)  \quad(1\le j\le s).
$$
Write $|x,y|_{i,j}$ for the image of $(x,y)$ in $\mathcal C_{i,j}$. 

For each $0\le i\le r$, let $\ell_b^i$ be the $b$-labelled line through $o_i=|0,0|_{i,0}$ in $\mathcal C_{i,0}$, and let $\ell_{c_0}^i$ be the $c_0$-labelled line through $o_i'=|\alpha_i,-\beta_i|_{i,0}$. We attach $\ell_b^i$ to the horizontal line $a=i$ in $\Pi_A$ and $\ell_{c_0}^i$ to the horizontal line through the origin in $\mathcal C_{i,1}$ identifying
$$
|n,0|_{i,0}\sim(n,i)_{\Pi_A},
\qquad
|\alpha_i,m-\beta_i|_{i,0}\sim|m,0|_{i,1}.
$$

For each $1\le j<s$, we attach the vertical $c_j$-labelled line through the origin in $\mathcal C_{i,j}$ to the horizontal $c_j$-labelled line through the origin in $\mathcal C_{i,j+1}$ by identifying $|0,m|_{i,j}$ and $|m,0|_{i,j+1}$.

Finally, we attach the vertical $c_s$-labelled line through the
origin in $\mathcal C_{i,s}$ to the vertical line $d=i$ in $\Pi_D$
by identifying $|0,m|_{i,s}$ and $(i,m)_{\Pi_D}$, where $n,m\in\Z$. 

Denote the resulting square complex by $Y_N$. Each gluing preserves
labels and orientations, so the covering maps from the pieces to the
two-generator tori define a map from $Y_N$ to the Salvetti complex
$X_N$ of $P_N$.

Just as in Section \ref{sec:construction-of-Y}, one can show that 
\begin{equation}\label{eq:long-cycle-incidence}
 (n,i)_{\Pi_A}=(i,m)_{\Pi_D}\text{ in }Y_N
 \quad\Longleftrightarrow\quad
 (n,m)=(\alpha_i+p_i k,\beta_i+q_i k)\text{ for some }k\in\mathbb Z.
\end{equation}
In particular, the auxiliary chain identifies $(n,0)_{\Pi_A}$ with
$(0,n)_{\Pi_D}$ for each $n\in\mathbb Z$. Write $v_n$ for this vertex
and $o=v_0$ for the basepoint.

As in Section~\ref{sec:local-isometry}, the map $Y_N\to X_N$ is a local isometry, $\pi_1(Y_N,o)$ embeds in $P_N$ and its universal cover embeds as a convex subcomplex of the universal cover of $X_N$.

\subsection{Free subgroup of \texorpdfstring{$G_N$}{G\_N} with undecidable membership}

Let $U_N$ be the image of $\pi_1(Y_N,o)$ in $P_N$. Let
$$
 h_i=b^{p_i}c_0^{-q_i},\qquad \Delta_i=b^{\alpha_i}c_0^{-\beta_i},
 \qquad g_i=a^i \Delta_i\qquad(0\le i\le r).
$$
A free basis of $U_N$ consists of the words
\begin{equation}\label{eq:long-cycle-U}
z_i=a^ih_i a^{-i}, \qquad  v_{i,j}=g_i(c_{j-1}c_j^{-1})g_i^{-1}, \qquad
 w_i=g_i d^{-i} \qquad (0\le i\le r,\ 1\le j\le s).
\end{equation}
Indeed, start with $\Pi_A$ and attach the cylinders of each chain in order. Each new cylinder meets the preceding union in its incoming line. The loop around the first cylinder is based at the end of $a^i$ and gives $z_i$. The path $g_i$ leads to the vertex $|\alpha_i,-\beta_i|_{i,0}$, which is identified with the origin of every subsequent cylinder of the chain. Going once around $\mathcal C_{i,j}$ from that vertex gives the based generator $v_{i,j}$ after returning along $g_i^{-1}$. Finally $\Pi_D$ is attached along $r+1$ disjoint lines. Using the line of the
auxiliary chain to join $\Pi_D$ to the preceding union gives one more free generator $w_i$ for each $1\le i\le r$ (note that $w_0=1$ by definition), exactly as in Section~\ref{sec:free-basis-of-u}. Consequently
$$
\rank(U_N)=(r+1)(s+1)+r=(r+1)(N-3)-1.
$$

Write $\lambda_n$ for the $a$-line in $\Pi_A$ with $b$-coordinate $n$, and $\mu_m$ for the $d$-line in $\Pi_D$ with $c_s$-coordinate $m$. These lines contain every $a$- and $d$-edge of $Y_N$. The auxiliary chain contains the vertex $v_n$ which is common to $\lambda_n$ and $\mu_n$. Equation \eqref{eq:long-cycle-incidence} says that the $i$-th chain gives a common vertex $(n,i)_{\Pi_A}=(i,m)_{\Pi_D}$ exactly when the $i$-th rule pairs $n$ with $m$. In particular, the path from $v_n$ to $v_m$ for one rule still has label $a^i d^{-i}$. The added cylinders identify the two endpoints of these segments at a single vertex.

Arguing as in Section \ref{sec:ad-paths-and-double-cosets}, we obtain

\begin{lemma}\label{lem:long-cycle-double-cosets}
The vertices $v_n$ and $v_m$ are joined by an $(a,d)$-path in $Y_N$
if and only if  $n\sim_{\calR}m$ if and only if $U_Nb^n L=U_Nb^m L$.
\end{lemma}

\begin{theorem}\label{thm:long-cycle-membership}
For every $N\ge6$, the raag $G_N$ contains a fixed finitely generated
free subgroup $H_N$ for which membership is undecidable for the
words $b^{-k}eb^k$, where $k\in\mathbb Z$.
\end{theorem}
\begin{proof}
Fix $\mathcal R,t$ as in Proposition~\ref{prop:arithmetic}, and set
$V_N=b^{-t}U_Nb^t$ and $H_N=\langle V_N,e\rangle$. 

Each $\mu_m$ contains the single auxiliary vertex $v_m$. Since the source arithmetic progressions are disjoint, by \eqref{eq:long-cycle-incidence}, each $\lambda_n$ contains at most one rule intersection, and $\lambda_t$ contains none. These are exactly the properties used in the proof of Proposition~\ref{prop:terminal-freeness}. The same argument shows that $V_N\cap L=\{1\}$. Therefore, $H_N\cong V_N*\langle e\rangle$. Its free basis consists of $e$ and the words $b^{-t}u b^t$ for $u$ in \eqref{eq:long-cycle-U}, so $\rank(H_N)=(r+1)(N-3)$.

Finally, by Proposition~\ref{prop:transfer} and Lemma \ref{lem:long-cycle-double-cosets} we have that 
$$
 b^{-k}eb^k\in H_N
 \quad\Longleftrightarrow\quad b^{-k}\in V_NL
 \quad\Longleftrightarrow\quad U_Nb^{t-k}L=U_Nb^tL
 \quad\Longleftrightarrow\quad t-k\sim_{\mathcal R}t.
$$
So, deciding $ b^{-k}eb^k\in H_N$  is equivalent to deciding membership in the connected component of $t$, which is undecidable.
\end{proof}

\section{Mapping class groups} \label{sec:mapping-class-groups}

In this section we prove Theorems \ref{thm:orientable-surface-membership} and \ref{thm:nonorientable-surface-membership}.

Suppose first that $\sigma(S)\le1$. The corresponding mapping class group is virtually free by the low-complexity computations in \cite[Chapter~2]{FarbMargalit}, hence has decidable membership problem.

Now suppose that $\sigma(S)\ge2$. Since the pentagon is isomorphic to its complement, by \cite[Lemma~3.27]{KatayamaKuno} there is an induced pentagon in the curve graph of $S$. By \cite[Proposition~7.16]{Koberda}, there is an embedding $\iota:G\hookrightarrow M$ sending the standard generators of $G$ to suitable powers of Dehn twists about the five curves. This finishes the proof of Theorem \ref{thm:orientable-surface-membership}

We now consider non-orientable case.  The groups $\Mod(N_{1,p})$ for $p\le3$, $\Mod(N_{2,p})$ for $p\le1$, and $\Mod(N_{3,0})$ are virtually free, see \cite[Sections~2 and~6]{HidberEtAl2022}. The subgroup membership problem is therefore decidable in these cases.

Consider $N_{2,2}$. The Birman exact sequence gives
$$
1\longrightarrow \pi_1(N_{2,1})
\longrightarrow \PMod(N_{2,2})
\longrightarrow \Mod(N_{2,1})
\longrightarrow 1,
$$
where $\PMod$ is the subgroup fixing each puncture. Here $\pi_1(N_{2,1})\cong F_2$, and $\Mod(N_{2,1})$ has an infinite cyclic subgroup of index four \cite[Theorem~A.5]{Stukow2006} generated by a Dehn twist. Its full preimage $Q$ has index eight in $M$ and fits into an exact sequence
$$
1\longrightarrow F_2\longrightarrow Q
\longrightarrow\Z \longrightarrow 1.
$$
This is a split extension and $Q\cong F_2\rtimes_\phi\mathbb Z$ for some $\phi\in\Aut(F_2)$. Every outer automorphism of $F_2$ is induced by a homeomorphism of a torus with one boundary component \cite[Section~3.1]{BridsonReid2020}. Thus $Q$ is the fundamental group of the compact $3$-manifold obtained as the mapping torus of such a homeomorphism. The subgroup membership problem in $Q$ is therefore decidable by \cite[Theorem~1]{FriedlWilton2016}. Since this decidability passes to finite extensions \cite[Lemma~12]{FriedlWilton2016}, it also holds in $M$.

In all the remainging cases we now construct an embedding of the pentagon raag using the blowup homomorphism of \cite[Proposition~2.4]{ParisSzepietowski2015}.
This homomorphism replaces selected punctures by crosscaps and is injective. We apply it to pure mapping classes, relative to the punctures that remain unchanged.

If $g=1$ or $g=2$, start with the orientable surface $S_{0,p+g}$ and replace $g$ punctures by crosscaps. If $g\ge 3$, start with $S_{1,p+g-2}$ and replace $g-2$ punctures by crosscaps. In either case the resulting surface is $N_{g,p}$. We have an embedding
$$
\PMod(S)
\hookrightarrow \Mod(N_{g,p}),
$$
where $S=S_{0,p+g}$ for $g=1,2$ and $S=S_{1,p+g-2}$ for $g\ge 3$.
Our assumptions imply that $p+g\ge 5$ in the case $g=1,2$, and $p+g-2\ge 2$ if $g\ge 3$. In both cases, the orientable surface $S$ satisfies $\sigma(S)\ge2$. Theorem~\ref{thm:orientable-surface-membership} therefore gives an embedding $G\hookrightarrow \Mod(S)$. As we argued above, this embedding maps each standard generator of $G$ to a power of a Dehn twist about a curve in $S$ is an embedding of $G$ into this mapping class group by powers of Dehn twists. These twists fix every puncture individually, so the image lies in $\PMod(S)$. Thus, $G\hookrightarrow\PMod(S)\hookrightarrow\Mod(N_{g,p})=M$.

\end{document}